\documentclass[11pt]{article}
\usepackage{amsfonts,latexsym,rawfonts,amsmath,amssymb,amsthm,graphicx}
\numberwithin{equation}{section}
\newtheorem{theorem}{Theorem}[section]
\newtheorem{lem}[theorem]{Lemma}
\newtheorem{thm}[theorem]{Theorem}

\def\endproof{$\hfill\Box$\\}

\def\div{{\rm div\,}}

\title{A Bernstein type theorem for graphic self-shrinkers with flat normal bundle}
\author{Yong Luo\footnote{The author is supported by the DFG Collaborative Research Center SFB/TR71.}}
\date{}
\begin{document}
\maketitle
\begin{abstract}In this note we will prove that an $n$ dimensional graphic self-shrinker in $R^{n+m}$ with flat normal bundle is a linear subspace. This result is a generalization of the corresponding result of Lu Wang in codimension one case.
\end{abstract}
\section{Introduction}
Let $\{M_t\}_{T>t>0}$ be a family of $n$ dimensional submanifolds in $R^{n+m}$ and let $H_t$ be the mean curvature vector of $M_t$ in $R^{n+m}$, then $\{M_t\}_{T>t>0}$ is said to be moving by mean curvature if
\begin{eqnarray}\label{mean curvature flow}
\frac{d}{dt}M_t=H_t.
\end{eqnarray}
Let $M$ be a $n$ dimensional submanifold in $R^{n+m}$ and $NM$ be its normal bundle. $M$ is said to the a self-shrinker if it satisfies
\begin{eqnarray}
H=-\frac{1}{2}\overrightarrow{x}^N,
\end{eqnarray}
where $\overrightarrow{x}$ is the position vector of $R^{n+m}$ and $\overrightarrow{x}^N$ is the normal component of $\overrightarrow{x}$ in $NM.$

Self-shrinkers are the simplest solutions to the mean curvature flow and they are important in the singularity analysis of the mean curvature flow, see for example \cite{CM2}, \cite{CM3}, \cite{Hui1} and \cite{Hui2}.

On the other hand, self-shrinkers can be viewed as minimal surfaces endowed with the the weighted metric $e^{-\frac{|\overrightarrow{x}|^2}{2n}}\delta_{ij}$, see \cite{Ang}, \cite{CM2} and \cite{CM3}. In the theory of minimal surfaces the Bernstein theorem for entire minimal graphs played a fundamental role (see, for example \cite{CM1} and \cite{Oss}), and so it is nature to ask whether there are Bernstein type theorems for graphic self-shrinkers. Ecker and Huisken proved that $n$ dimensional smooth self-shrinkers in $R^{n+1}$, which are entire graphs and have at most polynomial volume growth, are hyperplanes, in the appendix of \cite{EH}. The volume growth condition is removed by Wang \cite{W}. She proved that an $n$ dimensional self-shrinker of entire graph in $R^{n+1}$ has polynomial volume growth by using a calibration argument and a $L^\infty$ estimate to the graph function.

In this note, we will prove an Bernstein type theorem for graphic self-shrinkers with flat normal bundles of arbitrary codimension, which generalizes the result of Lu Wang. The main step is to get a stability inequality for the self-shrinker with the weighted measure $e^\frac{-|\overrightarrow{x}|^2}{4}d\mu$, by combining the ideas from \cite{xin} and \cite{W}. Our theorem is stated as follows.
\begin{thm}\label{main theorem}
Let $M$ be an $n$ dimensional graphic self-shrinker in $R^{n+m}$ with flat normal bundle, then $M$ is a linear subspace.
\end{thm}
Comparing to the proof of the Bernstein theorem for minimal graphs with flat normal bundle in \cite{SWX}, our proof of the above theorem is quite simple. The reason is that in our case we can get a "weighted stability inequality"(Lemma \ref{Stability inequality}) with weight $e^{-\frac{|\overrightarrow{x}|^2}{4}}$ for self-shrinkers, and then by the volume growth estimate of Ding and Xin (Theorem \ref{volume growth}) the right hand side of the weighted stability inequality tends to zero by choosing appropriate cut off functions.

For other Bernstein type theorems for graphic self-shrinkers, we refer to \cite{CCY}, \cite{DW}, \cite{DX2}, \cite{DXY} and \cite{HW}.

\textbf{Acknowledgement} I would like to thank my advisor, Prof. Guofang Wang for introducing me to this problem and discussions.
\section{Preliminaries}
\subsection{Basic equations from submanifold theory}
Let $M\hookrightarrow \overline{M}$ be an isometric immersion with the second fundamental form $B$, which can be viewed as a section of  the vector bundle $Hom(\odot^2TM, NM)$ over $M$, where $TM$ and $NM$ are the tangent bundle and normal bundle of $M$, respectively. Then the second fundamental form, the curvature tensor of the submanifolds, curvature tensor of the normal bundle and that of the ambient manifold, satisfy the Gauss equation, the Codazzi equation, and the Ricci equation, as follows.
\begin{eqnarray*}
\langle R_{XY}Z,W\rangle&=&\langle\overline{R}_{XY}Z,W\rangle-\langle B_{XW},B_{YZ}\rangle+\langle B_{XZ},B_{YW}\rangle,
\\(\nabla_XB)_{YZ}&-&(\nabla_YB)_{XZ}=-(\overline{R}_{XY}Z)^N,
\\ \langle R_{XY}\mu,\nu\rangle&=&\langle\overline{R}_{XY}\mu,\nu\rangle-\langle B_{Xe_i},\mu\rangle\langle B_{Ye_i},\nu\rangle-\langle B_{Xe_i},\nu\rangle\langle B_{Ye_i},\mu\rangle,
\end{eqnarray*}
where $\{e_i\}$ is a local orthonormal frame on $M$; $X, Y, Z\in TM$ and $\mu, \nu\in NM$.

In particular if $\overline{M}$ is the Euclidean space and $M$ has flat normal bundle, then the Ricci equation will be
\begin{eqnarray}\label{Ricci equation}
\langle B_{Xe_i},\mu\rangle\langle B_{Ye_i},\nu\rangle=\langle B_{Xe_i},\nu\rangle\langle B_{Ye_i},\mu\rangle.
\end{eqnarray}

\subsection{Volume growth of self-shrinkers}
The following theorem estimates the volume growth of self-shrinkers with any codimension (\cite{DX1}).
\begin{thm}[Ding-Xin]\label{volume growth}
Any complete non-compact properly immersed $n$ dimensional self-shrinker $M$ in $R^{n+m}$ has Euclidean volume growth at most. Precisely, $\int_{D_r}d\mu\leq Cr^n$ for $r\geq1$, where $C$ is a constant depending on $n$ and the volume of $D_{8n}$, and $D_r=M\cap B_r$, $B_r$ is the ball of $R^{n+m}$ centered at the origin with radius $r$.
\end{thm}

It is proved in \cite{CZ} that the inverse is also true. That is any complete self-shrinker in $R^{n+m}$ with Euclidean volume growth must be proper.

\section{Proof of theorem \ref{main theorem}}
For vectors $a_1,...,a_n$ and $b_1,...,b_n$ in $R^{n+m}$, we let $A=a_1\wedge...\wedge a_n$ and $B=b_1\wedge...\wedge b_n$ and define their inner product by
$$\langle A,B\rangle=\det(\langle a_i,b_j\rangle).$$
Let $M$ be an $n$ dimensional complete self-shrinker of $R^{n+m}$ with flat normal bundle. For a point $x\in M$ we chose an orthonormal frame field $\{e_i,e_\alpha\}$ such that $e_i\in TM$ and $e_\alpha\in NM$. Fix a $n$-vector $A=a_1\wedge...\wedge a_n$. We define a function on $M$ by
\begin{eqnarray}\label{The function}
\omega=\langle e_1\wedge...\wedge e_n,a_1\wedge...\wedge a_n\rangle=\det(\langle e_i,a_j\rangle).
\end{eqnarray}
Then we have
\begin{eqnarray*}
e_i(\omega)&=&\sum_j\langle e_1\wedge...\wedge D_{e_i}e_j\wedge...\wedge e_n,a_1\wedge...\wedge a_n\rangle
\\&=&\sum_j\langle e_1\wedge...\wedge (D_{e_i}e_j)^T\wedge...\wedge e_n,a_1\wedge...\wedge a_n\rangle
\\&+&\sum_j\langle e_1\wedge...\wedge (D_{e_i}e_j)^N\wedge...\wedge e_n,a_1\wedge...\wedge a_n\rangle
\\&=&\sum_{\alpha, j}h_{\alpha ij}\langle e_1\wedge...\wedge e_\alpha\wedge...\wedge e_n,a_1\wedge...\wedge a_n\rangle
\end{eqnarray*}
and so
\begin{eqnarray*}
\Delta \omega=-|B|^2\omega +\sum_{\alpha, i, j}h_{\alpha iij}\langle e_1\wedge...\wedge e_\alpha\wedge...\wedge e_n,a_1\wedge...\wedge a_n\rangle
\\+\sum_{\alpha, \beta, i, j, k}\langle e_1\wedge...\wedge h_{\alpha ij}e_\alpha\wedge...\wedge h_{\beta ik}e_\beta\wedge...\wedge e_n,a_1\wedge...\wedge a_n\rangle
\\=-|B|^2\omega +\sum_{\alpha,i,j}h_{\alpha iij}\langle e_1\wedge...\wedge e_\alpha\wedge...\wedge e_n,a_1\wedge...\wedge a_n\rangle
\\+\sum_{\alpha<\beta,i,j,k}(h_{\alpha ij}h_{\beta ik}-h_{\beta ij}h_{\alpha ik})\langle e_1\wedge...\wedge e_\alpha\wedge...\wedge e_\beta\wedge...\wedge e_n,a_1\wedge...\wedge a_n\rangle.
\end{eqnarray*}
Because $M$ is a submanifold in $R^{n+m}$ with flat normal bundle, by (\ref{Ricci equation}) we have
$$\sum_{\alpha<\beta,i,j,k}(h_{\alpha ij}h_{\beta ik}-h_{\beta ij}h_{\alpha ik})=0$$
and therefore we have
$$\Delta \omega=-|B|^2\omega +\sum_{\alpha,i,j}h_{\alpha iij}\langle e_1\wedge...\wedge e_\alpha\wedge...\wedge e_n,a_1\wedge...\wedge a_n\rangle.$$
On the other hand,
$$H=-\frac{1}{2}\overrightarrow{x}^N=-\frac{1}{2}\sum_\alpha\langle\overrightarrow{x},e_\alpha\rangle e_\alpha,$$
and so
\begin{eqnarray*}
\sum_{\alpha,i,j}h_{\alpha iij}=\sum_{\alpha,j}e_j\langle H,e_\alpha\rangle
=-\sum_{\alpha,j}\frac{1}{2}e_j\langle\overrightarrow{x},e_\alpha\rangle
=\sum_{\alpha,i,j}\frac{1}{2}h_{\alpha ji}\langle\overrightarrow{x},e_i\rangle,
\end{eqnarray*}
which implies that
\begin{eqnarray*}
\sum_{\alpha,i,j}h_{\alpha iij}\langle e_1\wedge...\wedge e_\alpha\wedge...\wedge e_n,a_1\wedge...\wedge a_n\rangle=\frac{1}{2}\sum_i\langle(e_i\omega) e_i,\overrightarrow{x}\rangle=\frac{1}{2}\langle \nabla\omega,\overrightarrow{x}\rangle.
\end{eqnarray*}
Finally we obtain the following formula for $\omega$
\begin{eqnarray}
\Delta\omega-\frac{1}{2}\langle\nabla\omega,\overrightarrow{x}\rangle+|B|^2\omega=0.
\end{eqnarray}
If $\omega>0$ on $M$, then let $g=\log\omega$ and $g$ satisfies the following equation
$$\Delta g+|\nabla g|^2-\frac{1}{2}\langle\nabla g,\overrightarrow{x}\rangle+|B|^2=0.$$
We have the following weighted stability inequality
\begin{lem}\label{Stability inequality}
Let $M$ be an $n$ dimensional complete self-shrinker in $R^{n+m}$ with flat normal bundle. If there is an $n$-vector $A$ such that the function $\omega$ defined by (\ref{The function}) is positive everywhere on $M$, then
$$\int_M|B|^2\eta^2e^{-\frac{|\overrightarrow{x}|^2}{4}}d\mu\leq\int_M|\nabla\eta|^2e^{-\frac{|\overrightarrow{x}|^2}{4}}d\mu,$$
where $\eta$ is any function with compact support on $M$.
\end{lem}
\proof Multiplying the equation of $g$ by $\eta^2e^{-\frac{|\overrightarrow{x}|^2}{4}}$ and integrating over $M$ gives
\begin{eqnarray*}
0&=&\int_M\eta^2\div(e^{-\frac{|\overrightarrow{x}|^2}{4}}\nabla g)+\int_M\eta^2(|\nabla g|^2+|B|^2)e^{-\frac{|\overrightarrow{x}|^2}{4}}
\\&=&-\int_M-2\eta\langle\nabla g,\nabla\eta\rangle e^{-\frac{|\overrightarrow{x}|^2}{4}}+\int_M\eta^2(|\nabla g|^2+|B|^2)e^{-\frac{|\overrightarrow{x}|^2}{4}}
\\&\geq&-\int_M(\eta^2|\nabla g|^2+|\nabla \eta|^2)e^{-\frac{|\overrightarrow{x}|^2}{4}}+\int_M\eta^2(|\nabla g|^2+|B|^2)e^{-\frac{|\overrightarrow{x}|^2}{4}}
\\&=&\int_M(-|\nabla \eta|^2+\eta^2|B|^2)e^{-\frac{|\overrightarrow{x}|^2}{4}},
\end{eqnarray*}
and the conclusion follows.
\endproof
 \textbf{Proof of theorem \ref{main theorem}:} Because $M$ is a graphic submanifold of $R^{n+m}$, we can find an $n$-vector $A$ such that $\omega$ is everywhere positive on $M$. Let $D_r=M\cap B_r$, where $B_r$ is the ball in $R^{n+m}$ centered at the origin with radius $r>1$, and we choose $0\leq\eta\leq1$ to be a function defined on $M$ which equals to $1$ on $D_r$ and equals to zero outside $D_{r+1}$, with first derivatives bounded by a constant $C$ independent of $r$. Then by lemma \ref{Stability inequality} we have
$$\int_{D_r}|B|^2e^{-\frac{|\overrightarrow{x}|^2}{4}}\leq\int_M\eta^2|B|^2e^{-\frac{|\overrightarrow{x}|^2}{4}}\leq\int_M|\nabla\eta|^2e^{-\frac{|\overrightarrow{x}|^2}{4}}\leq C\int_{D_{r+1}\setminus D_r}e^{-\frac{|\overrightarrow{x}|^2}{4}}.$$
Let $r\rightarrow\infty$ and by theorem \ref{volume growth} (note that $M$ is a graph in $R^{n+m}$, so it is proper) we have $B\equiv0$, which completes the proof of theorem \ref{main theorem}.
\
{}
\vspace{1cm}\sc

Yong Luo

Mathematisches Institut, Albert-Ludwigs-Universit\"at Freiburg,

Eckerstr. 1, 79104 Freiburg, Germany.

{\tt yong.luo@math.uni-freiburg.de}

\end{document}